\documentclass[11pt,a4paper]{article}
\usepackage{amsmath,amsfonts,amsbsy,latexsym,epsfig,amsthm,amssymb,graphics}

\title{Moduli Space of Cubic Surfaces as Ball Quotient via Hypergeometric Functions}

\author{Brent R. Doran}

\date{}

\newtheorem{Lem}{Lemma}

\newtheorem{Thm}{Theorem}
\newtheorem{Prop}{Proposition}
\newtheorem{Cor}{Corollary}
\newtheorem{Def}{Definition}
\newtheorem*{MainThm}{Theorem 3}

\theoremstyle{remark}
\newtheorem{remark}{Remark}

\newcommand{\ztwoz}{\ensuremath{\mathbb{Z}/2\mathbb{Z}}}
\newcommand{\Cmsm}{\ensuremath{\mathcal{M}^{\mbox{sm}}_{C,m}}}
\newcommand{\Csm}{\ensuremath{\mathcal{M}^{\mbox{sm}}_C}}

\else\oddsidemargin25mm\evensidemargin25mm\marginparwidth25mm\fi
\begin{document}

\maketitle

\abstract{We describe hypergeometric functions of Deligne-Mostow
type for open subsets of the configuration space of six points on
$\mathbb{P}^2$, induced from those for seven points on
$\mathbb{P}^1$. The seven point ball quotient example
$DM(2^5,1^2)$ does not appear on Mostow's original list \cite{M},
but does appear on Thurston's corrected version \cite{T}.  We show
that $DM(2^5,1^2)$ is a finite cover of the moduli space of cubic
surfaces $\mathcal{M}_C$ endowed with the ball quotient structure
$\Gamma_C \backslash \mathbb{B}^4$ of \cite{ACT}. This answers a
question of Allcock \cite{A} about the commensurability of
$\Gamma_C$ with the monodromy groups of Deligne-Mostow
hypergeometric functions. }

\section{Introduction}

A moduli space $\mathcal{M}$ of algebraic varieties can sometimes
be endowed with the structure of a Hermitian locally symmetric
space $\Gamma \backslash \mathcal{B}$ via a multi-valued period
map with discrete monodromy group $\Gamma$. In classical
terminology, the period map (or rather, its single-valued
$\Gamma$-invariant inverse) provides a ``uniformization'' of
$\mathcal{M}$ by the symmetric space $\mathcal{B}$. When
$\mathcal{B}$ is the complex $n$-dimensional hyperbolic space
form, i.e., the ``complex $n$-ball'' $\mathbb{B}^n$, then the
locally symmetric space is called a {\em ball quotient}.

Probably the best-known example is the uniformization of the
coarse moduli space of cubic (elliptic) curves in $\mathbb{P}^2$
by the upper half plane $\mathbb{H}^1 \cong \mathbb{B}^1$.  The
ball quotient structure on $\mathcal{M}$ is simply the $j$-line
($\mathbb{P}^1$ with two orbifold points --- of order 2 and 3 ---
and one puncture representing a cusp point). The $j$-function on
$\mathbb{H}^1$ is the inverse to the period map and is invariant
under the modular group $\Gamma = PSL_2(\mathbb{Z})$. Here the
period map is expressed in terms of classical hypergeometric
functions \cite{Y}.

One might hope for something similar with moduli of cubic surfaces
in $\mathbb{P}^3$. Unfortunately cubic surfaces have no non-zero
holomorphic 2-forms, and hence no non-trivial periods.

One option is to instead use the periods of naturally associated
varieties. The work of Allcock, Carlson, and Toledo \cite{ACT}
uses cubic threefolds that arise as tri-cyclic covers of
$\mathbb{P}^3$ branched over a cubic surface. They conclude that
the moduli space $\mathcal{M}_C$ of nodal cubic surfaces is a
$4$-ball quotient $\Gamma_C \backslash \mathbb{B}^4$, and that in
fact the GIT compactification (adding semi-stable points) of
$\mathcal{M}_C$ coincides with the Baily-Borel compactification
(adding cusp points) of the ball quotient.

Another option is to find a uniformization more directly using
hypergeometric functions. In the Deligne-Mostow framework, such
functions uniformize certain moduli spaces of $n$ points on
$\mathbb{P}^1$ by $\mathbb{B}^{n-3}$. Unfortunately, Allcock
\cite{A} showed that $\Gamma_C$ is not equal to any of the
Deligne-Mostow monodromy groups, so it is --- in that sense  --- a
new example and not of hypergeometric type.  However, his proof
leaves open whether $\Gamma_C$ is commensurable to a
Deligne-Mostow group.

\subsection{Motivation, intuition, and the argument in brief}

With this in mind, the motivating idea for how to relate
$\mathcal{M}_C$ and Deligne-Mostow hypergeometric functions is
straightforward.

$\mathcal{M}_C$ is $4$-dimensional.  The 4-dimensional
Deligne-Mostow ball quotients are moduli spaces of 7 points on
$\mathbb{P}^1$. All but three of these are descendants (i.e.,
realized as GIT-stable coincidences of subsets of 12 points, see
\cite{D1,D2} for details) of a single Deligne-Mostow 9-ball
quotient, which we call the Eisenstein ancestral example and
denote by $DM(1^{12})$. In particular, their monodromy groups are
defined (as automorphisms of lattices) over the Eisenstein ring
$\mathcal{E} = \mathbb{Z}[\omega], \omega^3 = 1$.

In \cite{ACT} it is shown that $\Gamma_C$ is likewise defined over
$\mathcal{E}$, and contains the monodromy group for a 6-point
Deligne-Mostow example $DM(2^6)$.  As the notation suggests,
$DM(2^6)$ is the descendant of $DM(1^{12})$ representing six pairs
of points sharing the same coordinates. However, only one 7-point
Eisenstein descendant contains $DM(2^6)$, namely $DM(2^5, 1^2)$.
(Observe that $DM(2^5,1^2)$ does not appear on the original list
of Mostow \cite{M}
--- descendants are perhaps the easiest way to see how the list misses
many examples --- although it does appear in Thurston's corrected
list \cite{T}.)  This strongly suggests that the ball quotient
monodromy group $DM_{\Gamma}(2^5,1^2)$ and $\Gamma_C$ are
naturally related. In Section \ref{Cubic} we give a geometric
proof that $DM_{\Gamma}(2^5,1^2)$ is a finite index subgroup of
$\Gamma_C$. The essential argument is summarized below.

A smooth cubic surface is isomorphic to $\mathbb{P}^2$ blown up at
6 generic points.  Label these points $m_i$.  The task is to
associate to these 6 points a set of 7 points on $\mathbb{P}^1$.
The first five $m_i$ uniquely determine an irreducible quadric
curve $Q$.  There are two tangent lines to $Q$ through $m_6$. Thus
$m_i, 1 \leq i \leq 5$ together with the two tangent points define
$7$ points on $Q$.  Projection from $m_6$ sends these to 7 points
on $\mathbb{P}^1$, where the first 5 are equally weighted (should
be considered as unordered points), as are the last 2 (the
ordering of the tangent lines is arbitrary).  These are the
symmetries of $DM(2^5,1^2)$.  The map is not quite invertible:
each fiber is a $(\mathbb{Z}/2\mathbb{Z})^4$ orbit.

Nodal cubic surfaces are isomorphic to $\mathbb{P}^2$ blown-up at
6 points on a common irreducible $Q$, with at most pairs of points
``colliding'' (being infinitely near). The limit in which $m_6$
lies on $Q$ thus corresponds precisely to the class of nodal cubic
surfaces, and projection now gives $6$ points on $\mathbb{P}^1$
with exactly the structure of $DM(2^6)$.

Nodal cubic surfaces are the only singular cubic surfaces
parametrized by $\mathcal{M}_C$.  So in fact, up to finite
identifications, $DM(2^5,1^2)$ encodes all the data of
$\mathcal{M}_C$.

The main result says more:

\begin{MainThm}
$DM(2^5,1^2)$ is a finite branched cover of $\mathcal{M}_C$.
Furthermore, the Deligne-Mostow ball quotient structure is the
same as the Allcock-Carlson-Toledo one induced by the covering
map. In particular, $DM_{\Gamma}(2^5,1^2) \subset \Gamma_C$ as a
finite index subgroup.
\end{MainThm}

\subsection{Contents of paper}

Section \ref{Eisenstein} summarizes the needed background material
from the Deligne-Mostow theory of hypergeometric functions.  It
focuses on a very important special example, from which many other
examples are derived, which we call the Eisenstein ancestral
example.  One can show, using the techniques of this section
applied to the generators that $DM_{\Gamma}(2^5,1^2)$ is a
subgroup of $\Gamma_C$, although it doesn't illuminate the
geometry (or show that it is of finite index).

Section \ref{Cubic} contains the proof of the main theorem. In
Section \ref{Remarks} our result is discussed in light of recent
work of Dolgachev, van Geeman, and Kondo [DGK].

The author noticed the correspondence between the Deligne-Mostow
4-ball quotient $DM(2^5,1^2)$ and the Allcock-Carlson-Toledo
presentation of $\mathcal{M}_C$ upon reading the preprint versions
of \cite{ACT} and \cite{HL}, and attempted to enfold it in a
larger program of realizing moduli spaces as subball quotients of
Deligne-Mostow type as part of his thesis work \cite{D1}. Although
similar to the Hurwitz spaces and moduli of inhomogeneous binary
forms in \cite{D1,D2}, this example doesn't fit naturally in that
context. Nevertheless it appears to be of independent interest.

The author was unaware of the work of Dolgachev, van Geeman, and
Kondo until the conference on Geometry and Topology of Quotients
at the University of Arizona in December 2002 \cite{Dolg}. Their
original motivation was to understand the ball quotient structure
on $\mathcal{M}_C$ via periods of K3 surfaces, and their beautiful
paper \cite{DGK} comprehensively details the isomorphisms among
the moduli spaces they consider. Some similar ideas, from the
point of view of hypergeometric differential equations, can also
be found in \cite{Y2}, although as the authors themselves point
out the analysis therein is incomplete.

\section{Eisenstein Ancestral Example}\label{Eisenstein}

This section presents necessary background on the hypergeometric
functions, monodromy groups, and lattices involved in the main
construction.

Deligne and Mostow \cite{DM,M} characterize which moduli spaces of
GIT-stable $n$-punctured projective lines are uniformized by a
complex ball via hypergeometric functions.  The list of examples
they produce, however, has many natural inclusions. One example in
particular, which we call the Eisenstein ancestral example,
contains almost all of the Deligne-Mostow examples for $n \geq 7$
as natural ``descendants'' (see \cite{D1,D2} for technical
details).

A GIT moduli space of $n$ distinct points on $\mathbb{P}^1$ is
specified by a choice of ample line bundle on $(\mathbb{P}^1)^n$
that linearizes the action of $SL_2(\mathbb{C})$. This is
equivalent data to a list of $n$ positive integers $\mu = (\mu_1,
\ldots, \mu_n)$. We denote the GIT-stable moduli space by
$DM(\mu)$.

\begin{remark} Each $DM(\mu)$ is a quasi-projective variety together with
an explicit embedding.  The projective completion is the GIT
compactification obtained by adding the ``semi-stable but not
stable" points (up to equivalence). All of the $DM(\mu)$ for fixed
$n$ are isomorphic on the open subset representing collections of
{\em distinct} points on $\mathbb{P}^1$.
\end{remark}

\begin{Def}
The Eisenstein ancestral example is $DM(1^{12})$, i.e., $DM(1,
\ldots, 1)$.
\end{Def}

\begin{remark}
Passing from ordered to unordered points by quotienting by
$S_{12}$ yields the moduli space of degree 12 binary forms.
\end{remark}

Hypergeometric functions on $DM(1^{12})$ are defined as follows.
Denote by $\mathbb{P}^1_S$ the $12$-punctured projective line
$\mathbb{P}^1 \setminus S, S = \{s_1, \ldots, s_{12} \}$. Consider
the multi-valued form
$$\omega_{HG}(s_1, \ldots, s_{12}) = \prod_i
(z-s_i)^{-\frac{1}{6}} dz$$ on $\mathbb{P}^1_S$ (assuming only for
convenience of notation that $s_i \neq \infty, \forall i$). This
is the same data as:  a rank one complex local system, $l_S
\rightarrow \mathbb{P}^1_S$, with monodromy about any of the $s_i$
given by multiplication by $e^{\frac{\pi \imath}{3}}$.  It is easy
to see that the dual local system is just the conjugate local
system $\overline{l_S}$ (characterized by local monodromies
$e^{\frac{-\pi \imath}{3}}$).

\begin{Prop}[after \cite{DM,M}]
The local system valued cohomology $H^1(\mathbb{P}^1_S, l_S)$ is
10 dimensional.  It admits an anti-Hermitian perfect pairing with
compactly supported cohomology, which is defined over the
Eisenstein integers $\mathcal{E} = \mathbb{Z}[\omega]$ (for
$\omega$ a cube root of unity).  The associated Hermitian form is
Lorentzian, i.e., has signature $(1,9)$.  We call the rank 10
Lorentzian lattice $\Lambda_{\mathcal{E}}$.  In the Hodge
decomposition $H^{1,0} \oplus H^{0,1}$, the holomorphic subspace
is positive definite, and the anti-holomorphic subspace is
negative definite, with respect to the Hermitian form.
\end{Prop}

\begin{remark}
Alternatively, one could use intersection cohomology to get a true
intersection pairing on $IH^1$ \cite{D1,D2}.
\end{remark}

Now let the points of $S$ vary on $\mathbb{P}^1$.  Deligne and
Mostow show that what one hopes to be true is in fact true, namely
that the $H^1(\mathbb{P}^1_S,l_S)$ are the fibers of a local
system $\mathcal{L}$ over the parameter space of $n$ distinct
points, $(\mathbb{P}^1)^n \setminus \Delta$, where $\Delta$ is the
union of the subdiagonals.  The resulting monodromy group
$DM_{\Gamma}(1^{12})$ must preserve $\Lambda_{\mathcal{E}}$, hence
$DM_{\Gamma}(1^{12}) \subset Aut(\Lambda_{E})$.

They also observe that projectivizing the fibers of $\mathcal{L}$
yields a {\em canonical} flat projective space bundle.
$H^{1,0}(\mathbb{P}^1_S,l_S)$ is one dimensional, so its
projective image is a point.  By abuse of notation, the projective
monodromy group we also denote by $DM_{\Gamma}(1^{12})$, and
henceforth take that to be the default meaning of the notation.
The natural extension of the group by the permutation group of the
weights, here $S_{12}$, is denoted by
$DM_{\Gamma,\Sigma}(1^{12})$.

Choose coordinates on a fiber of the bundle over an arbitrary base
point.  Use flatness to canonically extend coordinates over the
bundle, up to the projective monodromy action. Thus the
projectivized holomorphic component of the cohomology defines a
multi-valued map, with $DM_{\Gamma}(1^{12})$-monodromy, from the
parameter space of $12$ distinct points to $\mathbb{B}^9 \subset
\mathbb{P}^9$. It is easily seen that this map is invariant with
respect to the natural $SL_2(\mathbb{C})$-action on the parameter
space, and so descends to the moduli space $DM(1^{12})$. We call
this map the hypergeometric function $HG(1^{12})$ on $DM(1^{12})$.
They show it extends naturally to the boundary, which corresponds
to GIT stable collisions of the 12 points on $\mathbb{P}^1$.  This
boundary locus is stratified, indexed by partitions of the 12
points into subsets with no more than 5 elements.  The image of
$HG(1^{12})$, upon extending to the boundary, is all of
$\mathbb{B}^9$.

\begin{Thm}[after \cite{DM,M}] \label{Thm:EisensteinHG}
$HG(1^{12})$ defines a complex ball uniformization of
$DM(1^{12})$.  Specifically, there is an analytic isomorphism
$DM(1^{12}) \cong DM_{\Gamma}(1^{12}) \backslash \mathbb{B}^9$.
Furthermore, permutations of the coordinates by $S_{12}$ are
automorphisms of the lattice.
\end{Thm}

\begin{Cor}[\cite{DM,M}]
The GIT compactification of $DM(1^{12})$ is isomorphic to the
Baily-Borel compactification of $DM_{\Gamma}(1^{12}) \backslash
\mathbb{B}^9$.
\end{Cor}

The closure of a boundary stratum is again a Deligne-Mostow
example with $DM_{\Gamma}(\mu)$ defined over $\mathcal{E}$. Its
image under $HG(1^{12})$ is a $DM_{\Gamma}(1^{12})$-orbit of a
subball of $\mathbb{B}^9$. Indeed, the boundary divisor associated
to the collision of two points is given by, on the ball quotient
side, the Lorentzian sublattice orthogonal to a distinguished
vector (see \cite{D1} for details). The process is iterated for
further collisions.  A simple count shows:

\begin{Prop}
There are 6 Deligne-Mostow 4-ball quotient examples that are
descendants of $DM(1^{12})$.  As $DM(2^5, 1^2)$ and
$DM(3,2^3,1^3)$ are not on Mostow's original list \cite{M}, that
list is incomplete.
\end{Prop}
\begin{remark}
Mostow compiled his list using the correct theory, but made
mistakes in his by-hand computations.  Thurston later \cite{T}
corrected the computation with a computer search, and his list is
now believed to be complete (although I am unaware of a direct
proof).
\end{remark}

Note that each of the Deligne-Mostow monodromy groups is a
subgroup of the projective automorphism group of a lattice over
$\mathcal{E}$.  For 4-ball quotients it turns out there is only
one such lattice.

\begin{Thm}[\cite{A}]
Up to isometry, there is a unique Lorentzian rank 5 lattice over
$\mathcal{E}$, which we denote by $\Lambda_5$.
\end{Thm}

Denote Eisenstein 7 point Deligne-Mostow monodromy groups by
$DM_{\Gamma,\mathcal{E}}(\mu)$. The theorem tells us that
$\Gamma_C$ and all the $DM_{\Gamma,\mathcal{E}}(\mu)$ are
subgroups of the same $PAut(\Lambda_5)$.  It then follows from
Theorem \ref{Thm:EisensteinHG} that the extension
$DM_{\Sigma,\Gamma,\mathcal{E}}(\mu)$ of the Deligne-Mostow
monodromy group by $\Sigma$ (the permutations of the points of
equal weight) is likewise a subgroup of $PAut(\Lambda_5)$.

\section{Moduli Space of Marked Cubic surfaces and $DM(2^5,1^2)$}\label{Cubic}

Consider the parameter space of $6$ points on $\mathbb{P}^2$. Let
$\Delta$ denote the divisor representing non-generic
configurations, i.e., when all six points lie on a common quadric
curve or if any three are collinear (including multiplicities).  A
smooth cubic surface is the blow up of $6$ generic marked points
on $\mathbb{P}^2$. By the ``moduli space of smooth marked cubic
surfaces", $\Cmsm$, we mean the quotient of the parameter space of
generic configurations by automorphisms of $\mathbb{P}^2$. That
is:
$$\Cmsm \cong ((\mathbb{P}^2)^6
\setminus \Delta) / \! / SL_3(\mathbb{C})$$

This quotient is 4 dimensional. However, projectively inequivalent
markings can represent isomorphic cubic surfaces.  It is
well-known that the Weyl group of $E_6$, classically presented as
automorphisms of a cubic surface that permute its 27 lines, acts
on a marked $\mathbb{P}^2$ via birational Cremona transformations
and that these yield isomorphic cubic surfaces. Indeed:

\begin{Lem}[\cite{ACT}]
$\Cmsm / W(E_6) \cong \Csm$, where $\Csm$ is the moduli space of
smooth cubic surfaces. Furthermore, $W(E_6)$ acts as a subgroup of
$\Gamma_C$.
\end{Lem}
\label{Lem:W(E6)}

 Basic facts about the geometry of the
elementary (3-point) Cremona transforms will be needed.
\begin{Lem}\label{Lem:Cremona}
An elementary Cremona transform, defined by blowing
up three points on $\mathbb{P}^2$ and then blowing down the proper
transforms of the three lines that connect the vertices, has the
following properties, up to projective transformations:
\begin{enumerate}
\item It sends a line through an edge of the triangle to the
opposite vertex, and sends any line through a single vertex to
itself, albeit non-trivially.

\item It preserves any quadric that passes through precisely two
of the vertices.

\item Lines which don't meet any of the three vertices are sent to
quadrics passing through the three vertices.
\end{enumerate}
\end{Lem}
\begin{proof}
The first follows essentially by definition.  The last two are
easily seen in coordinates.  Pick coordinates so that $m_1 =
[1,0,0], m_2 = [0,1,0],$ and $m_3 = [0,0,1]$.  Then the Cremona
transformation is $\psi_{123}: [x,y,z] \mapsto [yz,xz,xy]$.  Any
nondegenerate quadric passing through two chosen points may be
presented, via projective automorphisms, as the Veronese quadric
$y^2-xz =0$ with the marked points $m_1$ and $m_2$ as above, and
with the third vertex sent to $m_3$.  The Veronese quadric in
$\psi_{123}$-transformed coordinates is $Y^2 - XZ = xz(xz-y^2)=0$,
so indeed the quadric is preserved.

A line that doesn't intersect $m_1, m_2,$ or $m_3$ may be taken to
be $x+y+z = 0$.  Its image set under $\psi_{123}$ is $yz+xz+xy=0$,
a non-degenerate quadric.  Also, $m_1,m_2,$ and $m_3$ are fixed
and lie on the image quadric. So all six $m_i$ are on the image
quadric.
\end{proof}

In the following Lemma we construct a map, $\phi^6_7$, from
configurations of 6 generic points on $\mathbb{P}^2$ to
configurations of 7 points on $\mathbb{P}^1$.   The two subsequent
Lemmas use $\phi^6_7$ in order to realize an open subset of
$DM(2^5,1^2)/(\mathbb{Z}/2\mathbb{Z})$ as a branched cover of
$\Csm$. This is the central idea governing the whole geometric
construction.

\begin{Lem}
There is a morphism
$$\phi^6_7: ((\mathbb{P}^2)^6 \setminus \Delta) / \! / SL_3(\mathbb{C}) \cong
\Cmsm \rightarrow DM_{7, \mathcal{E}}(\mu)/(\ztwoz)$$ given
geometrically by associating $7$ points $p_i$ on a projective line
to six generic ordered points $m_i$ on $\mathbb{P}^2$. If
$DM_{7,\mathcal{E}}(\mu)$ is taken to be $DM(2^5, 1^2)$, then the
morphism is naturally $S_5$-equivariant.  The $S_5$ acts as
automorphisms of the ball quotient structure on $DM(2^5,1^2)$.
\end{Lem}
\begin{proof}
Because the $m_i$ are generic (not a configuration represented in
$\Delta$), the first five points $\{ m_1, \ldots, m_5 \}$ lie on a
unique irreducible quadric $Q$ and $m_6$ is not on $Q$.  There are
precisely two tangent lines to $Q$ that pass through $m_6$. Label
the tangent points $t_1, t_2$, for a total of seven labelled
points on $Q$.  Since the tangent lines are unordered the
labelling of the $t_i$ is arbitrary.  More precisely, there are 5
ordered and 2 unordered points.

Projection from $m_6$ determines 7 marked points $p_i$ on the
image $\mathbb{P}^1$, by $m_i \mapsto p_i, 1 \leq i \leq 5$, and
$t_j \mapsto p_{j+5}$.  There is a $\mathbb{Z}/2\mathbb{Z}$ action
interchanging $p_6$ and $p_7$ since they are unordered. Observe
that any one of $\{p_i, 1 \leq i \leq 5\}$ may equal $p_6$ and any
one other may equal $p_7$, but no other multiplicities are
allowed.

The map is by construction compatible with the $SL_3$ and $SL_2$
actions: Any $SL_2$ action on the projective line of course may be
realized as the restriction of an action on $\mathbb{P}^2$.  To
see the reverse direction in coordinates, we may use $SL_3$ to
take $Q$ to be the Veronese embedding $V: y^2 = xz$ and $m_6 =
(0,1,0)$.  The associated isotropy subgroup acts as automorphisms
of the projective line.

The $S_5$-equivariance follows for the example $DM(2^5,1^2)$
because it permutes the first five coordinates on each side.  It
acts as automorphisms of $\Lambda_5$ by restriction from Theorem
\ref{Thm:EisensteinHG}.
\end{proof}

The map is not one-to-one.  Let $DM_{\ztwoz}(2^5,1^2)$ denote
$DM(2^5,1^2)/(\ztwoz)$, where \ztwoz\ acts by permuting the two
coordinates of weight 1.  (Deligne-Mostow descendants inherit the
action and so we will use similar notation for them.)

\begin{Lem}
$\phi^6_7$ descends to a dominant $S_5$-equivariant injective
morphism
$$\Cmsm /(\ztwoz)^4
\longrightarrow DM_{\ztwoz}(2^5,1^2),$$ with image set a union of
three strata: $DM^{\mbox{gen}}_{\ztwoz}(2^5,1^2),
DM^{\mbox{gen}}_{\ztwoz}(3,2^4,1)$, and
$DM^{\mbox{gen}}_{\ztwoz}(3^2, 2^3)$. It is an isomorphism onto
its image.
\end{Lem}
\begin{proof}
  In the Veronese presentation above, the pair $\{p_6, p_7\}$ is
canonically identified with the pair of tangent points
$\{[1,0,0],[0,0,1] \}$.  If the remaining $p_i$ are distinct from
$p_6$ and $p_7$, then $p_i = [1,0,x^2_i]$, $x_i \neq 0$, $1 \leq i
\leq 5$, and each one is the image of precisely two points on $Q$,
namely, $q_{i1} = [1,x_i,x_i^2]$ and $q_{i2} = [1, -x_i, x^2_i]$.
Thus any fiber of $\phi^6_7$ admits a transitive
$(\mathbb{Z}/2\mathbb{Z})^5$ action.  However, the automorphism of
$\mathbb{P}^2$ changing the sign of the $y$ coordinate obviously
preserves the quadric $Q$, the point of projection $[0,1,0]$, and
the $p_i$, while simultaneously switching $q_{i1}$ and $q_{i2}$
for all $i$.  So the effective transitive action on a generic
fiber (i.e., over the open subset
$DM^{\mbox{gen}}_{\ztwoz}(2^5,1^2)$) is $(\ztwoz)^5/(\ztwoz) \cong
(\ztwoz)^4$.

If one (or two) of the $p_i, 1 \leq i \leq 5$, call it $p_j$ (and
$p_k$), coincides with $p_6$ or $p_7$ (or both), then the
description is identical except that $q_{j1} = q_{j2}$ (and
$q_{k1}=q_{k2}$). Of course, $(\ztwoz)^4$ still acts transitively,
but not freely. The case with a single point of multiplicity 2 is
the stratum $DM^{\mbox{gen}}_{\ztwoz}(3,2^4,1)$ and the case with
two such points is the stratum
$DM^{\mbox{gen}}_{\ztwoz}(3^2,2^3)$.

Quotienting by the $(\ztwoz)^4$ action makes the inverse
well-defined, yielding the isomorphism onto the image. As before,
the $S_5$ action permutes the $m_i, 1 \leq i \leq 5$ on the one
hand and the $p_i, 1 \leq i \leq 5$ on the other.
\end{proof}

\begin{Lem}\label{Lem:ztwozCremona}
The groups $(\ztwoz)^4$ and $S_5$ act via Cremona transforms as
isomorphisms of cubic surfaces.
\end{Lem}
\begin{proof}
Let $m_1, \ldots, m_5$ define the nondegenerate quadric $Q$. Then
$m_6$ lies off $Q$.  The line $\overline{m_i m_6}$ intersects $Q$
in $m_i$ and $n_i$, with $m_i = n_i$ if and only if
$\overline{m_im_6}$ is a tangent line to $Q$.  Recall that
$(\ztwoz)^4$ acts via $(\ztwoz)^5$, the $i^{th}$ factor
interchanging $m_i$ with $n_i$, up to a $\ztwoz$ projective
transformation that switches {\em all} the $m_i$ and $n_i$ for $i
\neq 6$.  Thus it suffices to exhibit Cremona transformations, up
to projective equivalence, that switch any four $m_i$ with $n_i$
for $i \neq 6$. Without loss of generality, we show this for $i =
1, \ldots, 4$.

The Cremona transformation $\psi_{126}$, composed with a
projective change of coordinates whose notation we supress,
preserves $m_1$, $m_2$, and $m_6$ but switches $m_3$, $m_4$, and
$m_5$ with their respective $n_i$.  This follows from Lemma
\ref{Lem:Cremona}, because $Q$ passes through $m_1$ and $m_2$ and
hence is sent to itself, but also lines through $m_6$ are sent to
themselves non-trivially, so it must be the case that $m_i
\leftrightarrow n_i$ for $i = 3$, $4$, and $5$.  Then $\psi_{346}
\cdot \psi_{126}$ switches $m_i$ and $n_i$ for $i = 1, \ldots, 4$.

$S_5$ can be understood similarly, but it is easier to observe
that it (indeed $S_6$) acts by permuting the coordinates and so is
obviously an automorphism of a cubic surface and manifestly a
subgroup of the $W(E_6)$ action from Lemma \ref{Lem:W(E6)}.
\end{proof}

Let us simpify notation a bit by defining:
$$DM^{\mbox{sm}}_{\ztwoz}(2^5,1^2) = DM^{\mbox{gen}}_{\mathbb{Z}/2\mathbb{Z}}(2^5,1^2) \cup
DM^{\mbox{gen}}_{\ztwoz}(3,2^4,1) \cup
DM^{\mbox{gen}}_{\ztwoz}(3^2, 2^3)$$ It then follows from Lemma
\ref{Lem:ztwozCremona} that:

\begin{Cor}
$DM^{\mbox{sm}}_{\ztwoz}(2^5,1^2) \cong
\Cmsm/(\mathbb{Z}/2\mathbb{Z})^4$ is an intermediate branched
cover of $\Csm$, as are their respective $S_5$ quotients.

\end{Cor}

Let $U^{\mbox{sm}}$ and $V^{\mbox{sm}}$ denote the open subsets of
$\mathbb{B}^4$ such that
$$\Gamma_C \backslash U^{\mbox{sm}}
\cong \mathcal{M}_C \ \ \ \mbox{and} \ \ \ DM_{\Gamma}(2^5,1^2)
\backslash V^{\mbox{sm}} \cong DM^{\mbox{sm}}_{\ztwoz}(2^5,1^2).
$$
Our goal is to show there is a finite branched covering map
$DM_{\Gamma}(2^5,1^2) \backslash \mathbb{B}^4 \rightarrow \Gamma_C
\backslash \mathbb{B}^4$.  Thus we must extend the existing
covering map by adding the boundary loci (arrangements of 3-balls)
to $U$ and $V$ in $\mathbb{B}^4$.

We can try to extend $\phi^6_7$ to strata of $\Delta$ that match
with the boundary strata of $DM(2^5,1^2)$. The most important
case, which we denote by $\Delta^{\mbox{gen}}_Q$, is when all six
$m_i$ lie on an irreducible quadric $Q$. Denote by $\Delta_Q$ the
partial closure that allows points of multiplicity two.

\begin{Lem}
$\phi^6_7$ extends, as an $S_5$ equivariant map, over $\Delta_Q$.
It is an isomorphism onto its image $DM(2^6) \subset
DM_{\ztwoz}(2^5,1^2)$.
\end{Lem}
\begin{proof}
Projection from $m_6$ identifies 6 points on $\mathbb{P}^1$, with
$m_6$ identified with the point at $\infty$.  Thus
$\Delta^{\mbox{gen}}_Q$ is mapped to $DM^{\mbox{gen}}(2^6)$.
Pairwise collisions are precisely the stability condition for
$DM(2^6)$ (inherited from that of $DM(1^{12})$, described in
Section \ref{Eisenstein}). So the descendants of $DM(2^6)$ make up
the rest of the image of $\Delta_Q$.  The map is clearly
invertible.  The $S_5$ equivariance is reordering of
coorindinates, as before.
\end{proof}

It is well-known that the moduli space of singular cubic surfaces
with only nodal singularities, $\mathcal{M}^{\mbox{nod}}_C$, is
isomorphic to the variety $DM_{\Sigma}(2^6)$. What is more, the
3-ball quotient structure on $\mathcal{M}^{\mbox{nod}}_C$ induced
by restriction from $\mathcal{M}_C \cong \Gamma_C \backslash
\mathbb{B}^4$ coincides with the ball quotient structure of
$DM_{\Sigma}(2^6)$ \cite[Section 9]{ACT}. Let
$\mathbb{B}^3_{DM(2^6)}$ denote a 3-subball of $\mathbb{B}^4$
which covers $DM(2^6)$.  Then, in particular, $\Gamma_C$ acts on
$\mathbb{B}^3_{DM(2^6)}$ in such a way that the orbit is precisely
the arrangement of subballs which make up the complement of
$U^{\mbox{sm}}$ in $\mathbb{B}^4$. Then as an immediate
consequence, we have:

\begin{Cor}
The branched covering extends to $$(DM_{\Gamma}(2^5,1^2)
\backslash V)/(\ztwoz) \longrightarrow \Gamma_C \backslash
\mathbb{B}^4 \cong \mathcal{M}_C$$
\end{Cor}

To establish our goal, we simply need to show that the $\Gamma_C$
orbit of $\mathbb{B}^3_{DM(2^6)}$ in $\mathbb{B}^4$ is in fact the
complement of $V^{\mbox{sm}}$ as well.  The subgroup of $\Gamma_C$
that acts non-trivially on cubic surfaces is the factor group
associated to $\Cmsm$, namely $W(E_6)$ acting via Cremona
transformations.  Thus it suffices to show that all the remaining
strata of $DM_{\ztwoz}(2^5,1^2)$ can be obtained as images of
subsets of $DM(2^6)$ under Cremona transformations in $W(E_6)$.

We consider another extension of $\phi^6_7$ from $\Cmsm$ to handle
the other codimension 1 stratum of $DM_{\ztwoz}(2^5,1^2)$, namely
$DM^{\mbox{gen}}_{\ztwoz}(4,2^3,1^2)$.  On $\mathbb{P}^2$, this
stratum corresponds to the condition: $m_6$ is collinear with
$m_i$ and $m_j$, for some $i \neq j \in 1, \ldots, 5$.

\begin{Lem}
Cremona transforms identify $DM^{\mbox{gen}}(2^6)$ with
$DM^{\mbox{gen}}(4,2^3,1^2)$
\end{Lem}
\begin{proof}
Consider a point in $DM^{\mbox{gen}}(4,2^3,1^2)$.  It represents a
configuration where three points, of which $m_6$ is one, are
collinear.  The result now follows by using the Cremona
transformation from the proof of item 3 in Lemma
\ref{Lem:Cremona}.
\end{proof}

So, the two extensions of $\phi^6_7$ we have described are
identified by the Cremona transform. We now show that same Cremona
transform also identifies the boundary of these two divisors.  We
check this fact by taking limits in the respective extended
domains of $\phi^6_7$. Indeed, it is sufficient to check for the
one main case, $DM(4,2^4)$, as is clear from the poset of
descendant strata.

The interpretation via the extension of $\phi^6_7$ for $DM(4,2^4)
\subset DM(2^6)$ is straightforward: simply collide two points of
the six on the quadric $Q$.  On the other hand, for $DM(4,2^4)
\subset DM_{\ztwoz}(4,2^3,1^2)$ one takes the other limit in which
there is only one ``tangent" point to $Q$, namely when it
degenerates into two lines.  The weight $4$ point on
$\mathbb{P}^1$ corresponds to either (I) two distinct points $m_i$
and $m_j$ ($i,j \neq 6$) collinear with $m_6$, or (II) two
overlapping points (the cases are equivalent up to the
$(\ztwoz)^4$ action). In each case, the Cremona transform defined
by the non-collinear triple results in a quadric with 5 marked
points: (I) One of $m_i$ or $m_j$ must also be collinear with a
side of the Cremona triangle (i.e., two elements of the triple),
and hence is identified with the opposite vertex for a total of
five points, one of double weight.  Also, all five points must lie
on the quadric that is the image of the line defined by $m_6$,
$m_i$ and $m_j$. (II) The line connecting $m_6$ and the doubled
point is sent to a quadric passing through the three vertices of
the Cremona triangle.  This proves:

\begin{Thm}\label{Thm:MainThm}
$DM(2^5,1^2)$ is a finite branched cover of $\mathcal{M}_C$.
Furthermore, the Deligne-Mostow ball quotient structure is the
same as the Allcock-Carlson-Toledo one induced by the covering
map. In particular, $DM_{\Gamma}(2^5,1^2) \subset \Gamma_C$ as a
finite index subgroup.
\end{Thm}
\begin{remark}
Allcock shows $\Gamma_C$ is unequal to any $DM_{\Gamma}(\mu)$, and
hence that $\Gamma_C$ is in some sense a ``new" example. Theorem
\ref{Thm:MainThm} answers his question \cite[Section 8]{A}
regarding commensurability.
\end{remark}

\section{Further Remarks}\label{Remarks}

Note that although we used $\phi^6_7$ and two of its extensions to
construct associated cubic surfaces and hence the finite cover, a
simultaneous extension in $(\mathbb{P}^2)^6 \setminus \Delta$ is
not well-defined. For example, the configurations with a point of
multiplicity 2 are interpretable in multiple ways in
$DM(2^5,1^6)$.

To correctly extend $\phi^6_7$, one needs to extend the domain of
definition by ``separating" these different interpretations.  They
are indexed by where the collision ``came from", be it quadric,
collinear, or a generic collision.

More formally then, the true extension of $\phi^6_7$ is to the
moduli space of marked nodal cubic surfaces compactified by 40
points \cite{N}. There are a total of 36 boundary divisors,
interchanged by the action of $(\ztwoz)^4 \ltimes S_5 \subset
W(E_6)$. The 36 are indexed as follows: (a) 1 in which the points
all lie on a quadric, (b) 10 such that $m_6$ is collinear with
some $m_i$ and $m_j$, (c) 10 such that three of the $m_k, k \neq
6$ are collinear, and (d) 15 that represent the ``six choose
three" pairwise collisions. If the full group $W(E_6)$ is allowed
to act, then all the divisors are identified, and the quotient is
$\mathcal{M}_C$.

The $S_5$ action identifies the 10 divisors of (b) with a single
locus corresponding to $DM(4,2^3,1^2)$.  It identifies those of
(c) with a single locus corresponding to $DM(2^6)$ (now thought of
as parametrizing the projections from $m_6$ of five points on a
degenerate quadric).

Rather than write out the full details, for a more formal
discussion  (one originally motivated by periods of K3 surfaces),
one should see the geometric proof in \cite{DGK}.  The authors
give a beautiful interpretation of this cover as the moduli space
of cubic surfaces with a marked line.  The $(\ztwoz)^4 \ltimes
S_5$ group action is interpreted as $W(D_5) \subset W(E_6)$.

Dolgachev, van Geeman, and Kond\={o} don't use the locus of 6
points on a quadric in constructing their isomorphism.   They use,
in effect, the alternate limit where $Q$ degenerates to a pair of
lines while $m_6$ remains off $Q$. The singular point is the
confluence of the tangent points, and so the image under
projection to $\mathbb{P}^1$ from $m_6$ has weight 2, yielding
$DM(2^6)$. This is the basis for their list of 19 strata on pp.
15-16.

\end{document}